\documentclass{article}[12pt]

\usepackage{amsmath}   
\usepackage{amssymb}
\usepackage{amsthm}
\usepackage{xcolor}

\begin{document}
	
\title{Elliptic non-degeneracy and spectral rigidity of classical billiards}

\author{Georgi Popov and Peter Topalov}

\date{}

\maketitle
	
\begin{abstract}
In this paper we show that the billiard ball map of the Liouville billiard tables of classical type on 
the ellipsoid is non-degenerate at the elliptic fixed point. As a corollary we obtain a spectral rigidity result.
\end{abstract}

\newtheorem{Th}{Theorem}  
\newtheorem{Prop}{Proposition}  
\newtheorem{Def}{Definition}  
\newtheorem{Lem}{Lemma}  
\newtheorem{Coro}{Corollary}
\newtheorem{Rem}{Remark}
  
\newcommand{\eqdef}{\stackrel{\rm def}{=}}  
\newcommand{\al}{\alpha}  
\newcommand{\be}{\beta}  
\newcommand{\ga}{\gamma}  
\newcommand{\Ga}{\Gamma}  
\newcommand{\ep}{\epsilon}  
\newcommand{\de}{\delta}  
\newcommand{\la}{\lambda}  
\newcommand{\te}{\theta}  
\newcommand{\om}{\omega}  
\newcommand{\si}{\sigma}  
\newcommand{\C}{{\mathbb C}}
\newcommand{\R}{{\mathbb R}}
\newcommand{\Z}{{\mathbb Z}}
\newcommand{\N}{{\mathbb N}}
\newcommand{\I}{{\mathcal I}}

\section{Introduction}
This article is a part of a project (cf. \cite{PT1}-\cite{PT5}) investigating the relationship between the dynamics of 
completely integrable or ``close'' to completely integrable billiard tables, the integral geometry on them, and the spectrum
of the corresponding Laplace-Beltrami operators. It is concerned with the spectral rigidity of the Laplace-Beltrami operator
associated with $C^1$-smooth (with respect to the parameter $t$) deformations  $X_t$, $t\in [0,1]$, of a smooth
billiard table $X_0$ in an ambient Riemannian  manifold $(\widetilde X,g)$ of dimension ${\rm dim}\, \widetilde X =2$
equipped with a Riemannian metric $g$, as defined in Sect. \ref{sect:rigidity-of-BT}. By definition, a billiard table
$X$ in $(\widetilde X,g)$ is a  smooth   compact submanifold of  $\widetilde X$  of dimension two with boundary
$\Gamma=\partial X$.  The elastic reflection of geodesics of  $(X,g)$ at $\Gamma$ determines
continuous curves on $X$ called {\em billiard trajectories}  as well as a discontinuous dynamical system on
the corresponding co-sphere bundle $S^\ast X$ -- the ``{\em billiard flow}'' consisting of broken bicharacteristics of
the Hamiltonian $H$ associated to $g$  via the Legendre transform. The latter induces  a discrete dynamical system $B$
called billiard ball map which is well-defined  and smooth on  the open  coball bundle ${\bf B}^\ast \Gamma$ of 
$\Gamma$  provided $\Gamma$ is locally strictly geodesically convex. In this case 
the map $B : {\bf B}^\ast \Gamma\to {\bf B}^\ast \Gamma$  is  exact symplectic. 
We are interested in the dynamical properties and the spectral rigidity of integrable and ``close'' to  integrable billiard tables
on the ellipsoid. 

\subsection{Liouville billiard tables}
The Liouville billiard tables  are the most important examples of integrable  billiard tables.
By definition (see \cite{PT1}), $(X,g)$ is a  {\em Liouville billiard table} if there exists a double covering map with two branched points,
\begin{equation}\label{eq:covering}
\tau : C\to X,
\end{equation}
where $C$ denotes the cylinder ${(\R/\Z)}\times [-N, N]$, $N>0$, coordinatized by the variables $x$ and $y$
respectively, so that the metric $\tau^*(g)$ and the integral $\tau^*(\I)$ have the following form on $C$,
\begin{eqnarray}\label{eq:g-Liouville_form}
dg^2&=&\big(f(x)-q(y)\big)(dx^2+dy^2)\\
d\I^2&=&\alpha\,dF^2+\beta\,dg^2\nonumber
\end{eqnarray}
where $\alpha\ne 0$ and $\beta$ are real constants and
\begin{equation}\label{eq:the_integral}
dF^2:=\big(f(x)-q(y)\big)\big(q(y)\,dx^2+f(x)\,dy^2\big)\,.
\end{equation}
In other words, the integral $d\I^2$ belongs to the pencil of $dg^2$ and $dF^2$.
Here $f\in C^\infty(\R)$ is 1-periodic, $q\in C^\infty([-N,N])$, and
\begin{itemize}  
\item[(i)] $f$ is  even, $f>0$ if $x\notin\frac{1}{2}{\Z}$, and        
$f(0)=f(1/2)=0$;  
\item[(ii)] $q$ is even, $q<0$ if $y\ne 0$,  $q(0)=0$ and $q^{''}(0)<0$;  
\item[(iii)] $f^{(2k)}(l/2)=(-1)^kq^{(2k)}(0)$,  $l=0,1$,    
for every natural $k\in{\N}$.  
\end{itemize}  
In particular, if  $f(x)\sim \sum_{k=1}^{\infty}\ f_kx^{2k}$ is the Taylor expansion of $f(x)$  
at $0$, then, by  (iii),  the Taylor expansion of $q(y)$ at $0$  
is $q(y)\sim \sum_{k=1}^{\infty}\ (-1)^k f_kx^{2k}$.  Moreover, $f^{''}(0)=-q^{''}(0)>0$. 

\begin{Rem}
In coordinates, the covering map $\tau : C\to X$ satisfies $\tau(x,y)=\tau(-x,-y)$.  
In this way, the branched points of $\tau$ are $F_1:=\tau(0,0)$ and $F_2:=\tau(1/2,0)$.
The metric \eqref{eq:g-Liouville_form} and the integral \eqref{eq:the_integral} on $C$ vanish at these points. 
All other points on $X$ are regular values of $\tau$ and the preimage of any regular value contains exactly  two
points.
\end{Rem}

\noindent In addition, we will assume that 
\begin{itemize}
\item[(iv)] $f$ has a symmetry with respect to the line $x=1/4$, i.e., $f(x)=f(1/2-x)$ for $x\in[0,1/2]$;
\item[(v)] $f$ has a {\em Morse singularity} at $x=1/4$ which amounts to $f''(1/4)<0$;
\item[(vi)]  the boundary $\Gamma$ is locally strictly geodesically convex, which amounts to $(\partial q/\partial y)(N)<0$. 
\end{itemize}
\noindent Liouville billiard tables satisfying the additional conditions (iv) - (vi) are said to be of {\em classical type}. 
 
 Note that then the line  (taken twice) on the cylinder $C$ corresponding to $x=1/4$ is a 
{\em  closed broken geodesic} $\gamma_1$ of $(X,g)$ with two vertices. In what follows we will study 
the (squared) billiard map
\begin{equation}\label{eq:B^2}
P:=B^2 : {\bf B}^*\Gamma\to{\bf B}^*\Gamma,\quad\Gamma:=\partial X,
\end{equation}
in a neighborhood of the fixed point corresponding to this particular broken geodesic. 
Here ${\bf B}^*\Gamma$ denotes the unit co-tangent bundle of the boundary of the billiard $(X,g)$.

According to \cite{PT1},  Proposition 2.9,   any Liouville billiard table  has the {\em string property} which 
means that any broken geodesic starting from the singular point $F_1\, (F_2)$ passes through 
$F_2\, (F_1)$ after the first reflection at the boundary and the sum of distances from any
point of $\Gamma$ to $F_1$ and $F_2$ is constant. Thus the Liouville billiard table can be 
regarded as a natural generalization of the billiard in the interior of the ellipse. 

\begin{Rem}\label{rem:involutions}
Consider the involutions $\sigma_j:C\to C$, $j=1,2$, given by $\sigma_1(x,y):= (1/2 - x,y)$ and
$\sigma_2(x,y):= (x,-y)$ for $(x,y)\in C$. The involutions $\sigma_1$ and $\sigma_2$ commute with each other
as well as with the map $(x,y)\mapsto (-x,-y)$. Moreover, the metric and the integral given by \eqref{eq:g-Liouville_form}
are invariant with respect to $\sigma_j$, $j=1,2$. Thus the involutions  $\sigma_j$, $j=1,2$, factor trough $\tau$ to 
smooth involutions ${\mathcal J}_{j}:X\to X$ which preserve the metric $g$ and the integral $\I$. 
Hence, for any  Liouville billiard table $(X,g)$ of classical type,  there is a group  ${\rm Iso}(X)\cong{\Z}_2\oplus{\Z}_2$ 
acting on $(X,g)$ by isometries, which is generated by ${\mathcal J}_{j}$, $j=1,2$. The set of fixed points of
${\mathcal J}_{1}$ coincides with the image of the  geodesic $\gamma_1$. The point $\varrho_0:=\tau(1/4,0)$ is a fixed
point of both ${\mathcal J}_{1}$ and ${\mathcal J}_{2}$. 
\end{Rem}

\noindent 
We are going to investigate the billiard ball map in a neighborhood of the vertices of  the bouncing ball  geodesic $\gamma_1$. 
Let 
\begin{equation}\label{eq:f_taylor}
f(x)=\alpha_0+\alpha_1(x-x_0)^2+\alpha_2(x-x_0)^4+O\big((x-x_0)^6\big)
\end{equation}
where $\alpha_0>0$ and $\alpha_1<0$ be the Taylor's expansion of $f$ at $x_0=1/4$.
Consider the fixed point  $p_*$ of the map \eqref{eq:B^2} that corresponds to the vertex $\tau(1/4,-N)\in\Gamma$ of
the bouncing ball geodesic $\gamma_1$ (the other fixed point $B(p_*)$ corresponds to the vertex $\tau(1/4,N)$). 
Let $\{(\theta,I)\}$ be action-angle variables of the completely integrable map $B^2$ in an open neighborhood of
$p_*$ in ${\bf B}^*\Gamma$ normalized by $I(p_*)=0$. 
In these coordinates, the map $P=B^2$ has the form 
\begin{equation}\label{eq:L}
(\theta,I)\mapsto\Big(\theta+\frac{dL}{dI}(I),I\Big)
\end{equation}
where the Hamiltonian $L(I)$ is chosen so that $L(0)=0$.  
The coordinates $\{(\theta,I)\in \R/\Z\times \R_+\}$, $\R_+:=[0,+\infty)$, can be considered as polar 
symplectic coordinates in a neighborhood of $p_*$. We have the following

\begin{Th}\label{th:hessian_elliptic_point}
For any Liouville billiard table of classical type the Hamiltonian $L$ is well defined 
and $C^\infty$-smooth in a neighborhood of zero in $\R_+$ and satisfies
\begin{equation}\label{eq:K'}
\frac{dL}{dI}(0)=-\frac{\sqrt{-\alpha_1}}{\pi}\int_{-N}^N\frac{dy}{\sqrt{\alpha_0-q(y)}}\neq 0
\end{equation}
and
\begin{equation}\label{eq:K''}
\frac{d^2L}{dI^2}(0)=\frac{\alpha_1}{4\pi^2}\left(2\int_{-N}^N\frac{dy}{\big(\alpha_0-q(y)\big)^{3/2}}-
\frac{3\alpha_2}{\alpha_1^2}\int_{-N}^N\frac{dy}{\sqrt{\alpha_0-q(y)}}\right).
\end{equation}
\end{Th}

Theorem \ref{th:hessian_elliptic_point} will be proven in Sect. \ref{sect:proof_hessian}. 

We are going to apply Theorem \ref{th:hessian_elliptic_point} in the case of Liouville billiard tables on the ellipsoid. 

\subsection{Liouville billiard tables on the ellipsoid}
Liouville billiard tables of classical type appear on the ellipsoid and on the surfaces of 
constant curvature (see \cite[\S 3]{PT1}). Let us consider the case of the ellipsoid
$E(a_0,a_1,a_2)$ given by the equation
\begin{equation}\label{eq:ellipsoid_introduction}
\frac{x_0^2}{a_0}+\frac{x_1^2}{a_1}+\frac{x_2^2}{a_2}=1, 
\end{equation}
where $\{(x_0,x_1,x_2)\}$ are the coordinates in $\R^3$ and $0<a_0<a_1<a_2$.
There are two types of Liouville billiard tables on \eqref{eq:ellipsoid_introduction}:
The billiards of the {\em first type} are given by the intersection of \eqref{eq:ellipsoid_introduction}
with the interior of the hyperboloid of one sheet
\begin{equation}\label{eq:confocal_family}
\frac{x_0^2}{a_0-\lambda}+\frac{x_1^2}{a_1-\lambda}+\frac{x_2^2}{a_2-\lambda}\le 1
\end{equation}
where $\lambda\in(a_0,a_1)$. Note that for any given value of the parameter $\lambda\in(a_0,a_1)$ 
there are two isometric billiards of this type -- one with $x_0>0$ and another one with $x_0<0$. 
The bouncing ball  geodesic $\gamma_1$ corresponds to the intersection of the billiard table with
the plane $\{x_2=0\}$.
The billiards of the {\em second type} are given by the intersection of \eqref{eq:ellipsoid_introduction} with 
the interior of the hyperboloid of two sheets \eqref{eq:confocal_family} where $\lambda\in(a_1,a_2)$.
As in the case of the billiards of the first type, for any given value of the parameter $\lambda\in(a_1,a_2)$
there are two isometric billiards of the second type -- one with $x_2>0$ and another one with $x_2<0$.
The  bouncing ball  geodesic $\gamma_1$ corresponds to the intersection of the billiard with
the plane $\{x_0=0\}$. We denote by $X^\lambda$, $\lambda\in(a_0,a_1)\sqcup(a_1,a_2)$, 
the family of such billiard tables equipped with the induced Euclidean metric from $\R^3$ and by
$L(I,\lambda)$ the corresponding Hamiltonian \eqref{eq:L}. We have the following

\begin{Th}\label{th:non-degeneracy}
Consider the billiard table $X^\lambda$ on the ellipsoid \eqref{eq:ellipsoid_introduction}. Then,
\begin{itemize}
\item[(i)] For any value of the parameter $\lambda\in(a_0,a_1)$ we   have 
$\frac{d^2L}{dI^2}(0, \lambda)<0$ and 
\[
0<-\frac{dL}{dI}(0, \lambda)<\sqrt{\frac{a_1}{a_2}}<1
\]
where $\lambda \mapsto -\frac{dL}{dI}(0, \lambda)$ is a strictly decreasing function of $\lambda\in(a_0,a_1)$ such that
$\frac{dL}{dI}(0, a_1-):=\lim_{a\to a_1, a<a_1}\frac{dL}{dI}(0, a)=0$;
\item[(ii)] For any value of the parameter $\lambda\in(a_1,a_2)$ we have  have
$\frac{d^2L}{dI^2}(0, \lambda)>0$ and 
\begin{equation*}
0<-\frac{dL}{dI}(0, \lambda)<\sqrt{\frac{a_2(a_2-a_0)}{a_0(a_1-a_0)}}
\end{equation*}
where $\lambda \mapsto-\frac{dL}{dI}(0, \lambda)$ is a strictly increasing function of $\lambda\in(a_1,a_2)$
such that $\frac{dL}{dI}(0, a_1+):=\lim_{a\to a_1,a>a_1}\frac{dL}{dI}(0, a)=0$ and
\begin{equation}\label{eq:right_bound2}
\sqrt{\frac{a_1(a_1-a_0)}{a_0(a_2-a_0)}}<-\frac{dL}{dI}(0, a_2-)\,.
\end{equation}
\end{itemize}
\end{Th}
The proof of Theorem \ref{th:non-degeneracy} will be given  in Sect. \ref{sect. LBT_on_ellipsoid}.

The following Corollary follows from the fact that $\frac{d^2L}{dI^2}(0)$ is
strictly negative for the billiards of the first type and strictly positive for the billiards of 
the second type.

\begin{Coro}
There is no Liouville billiard table of the first type that is isometric to a Liouville billiard table of the second type.
\end{Coro}

 By definition, the  fixed point $p_\ast$ of the map $P=B^2$ corresponding to a Liouville billiard table
$X^\lambda$   is {\em non-degenerate} if $\frac{d^2L}{dI^2}(0,\lambda)\ne 0$. In this case we  say as well that the map $P$ is
{\em twisted} at $p_\ast$. This condition allows one to apply the Kolmogorov-Arnold-Moser (KAM)  theorem for small perturbations
of the symplectic map $P$. The  fixed point is {\em elliptic} if $\frac{dL}{dI}(0,\lambda)\notin \Z$ and {\em $4$-elementary} if
in addition 
\begin{equation}\label{eq:4-elementary}
\frac{dL}{dI}(0,\lambda)\!\!\!\mod 1\notin\left\{\frac{1}{4},\frac{1}{3},\frac{1}{2},\frac{2}{3},\frac{3}{4}\right\}
\end{equation}
where $\frac{dL}{dI}(0,\lambda)\!\!\!\mod 1$ denotes the fractional part of $\frac{dL}{dI}(0,\lambda)$. 

We have the following 

\begin{Coro}\label{coro:non-degeneracy}
The  fixed point $p_\ast$ of the map $P=B^2$ corresponding to any Liouville billiard table  on the ellipsoid is always  non-degenerate.
If the axes of the ellipsoid satisfy $\sqrt{\frac{a_1}{a_2}}<1/4$ then any Liouville billiard table $X^\lambda$ 
of the first type is elliptic and $4$-elementary. More generally,  the fixed point $p_\ast$ of a Liouville billiard table 
of first  type is always elliptic and it is $4$-elementary for all but no more than five  billiards. The  fixed point  $p_\ast$  of a Liouville billiard table 
of the second  type is elliptic and $4$-elementary for all but no more than finitely many  billiards. 
\end{Coro}

The five exceptional billiards of first type appearing in Corollary \ref{coro:non-degeneracy} correspond to 
the values of $\lambda\in(a_0,a_1)$ for which \eqref{eq:4-elementary} does not hold. 
Note that in this case we have  $-\frac{dL}{dI}(0,\lambda)\in(0,1)$ by Theorem \ref{th:non-degeneracy} (i).
In contrast, in the case of the billiards of the second type, the estimate \eqref{eq:right_bound2} implies that 
$-\frac{dL}{dI}(0,a_2-)$ can be made greater than any a priori given positive number. 
In order to see this we set $a_0:=1$, $a_2:=a_1+1$, and then take $a_1\to\infty$ in the expression on 
the left side of \eqref{eq:right_bound2}. This implies that the number of the values of $\lambda\in(a_1,a_2)$ 
for which \eqref{eq:4-elementary} does not hold is finite but it could be arbitrarily large (and not just five as in the case of
 the billiards of the first type).
\subsection{Applications to the spectral rigidity} The main result about the spectral rigidity of Liouville billiard tables 
on the ellipsoid as well as an overview of preceding results  is obtained in Sect. \ref{sect:rigidity-of-BT}. 
We will give a brief description of it. 

Let us denote by ${\mathcal E}$ the family of all Liouville billiard tables  of first or second type on the ellipsoid
$E(a_0,a_1,a_2)$ exept the finitely many billiard tables described in Corollary \ref{coro:non-degeneracy}.
The corresponding involutions ${\mathcal J}_{j}$, $j=1;2$, depend only on type of the billiard table
(see \eqref{eq:involutions_ellipsoid}). 
We denote by  ${\mathcal B}$ the class of all billiard tables $X$ on the ellipsoid $E(a_0,a_1,a_2)$ which are
invariant with respect to the involutions ${\mathcal J}_{j}$, $j=1;2$. By definition, for any billiard table
$X\in {\mathcal E}$ the elliptic geodesic $\gamma_1$ is $4$-elementary and the corresponding
Birkhoff normal form of $P=B^2$ is non-degenerate. This is an open condition in the 
$C^5$ topology, thus there exists a neighborhood ${\mathcal O}(X)$ of $X$ in  ${\mathcal B}$ in the 
the corresponding $C^5$ topology, such that any billiard $Y\in {\mathcal O}(X)$ enjoys that property.  

We consider $C^1$ families  $X_t$, $t\in [0,1 ]$, of billiard tables on the ellipsoid (see Definition \ref{def:C-1-families})
and denote  by $\Delta_t$ the corresponding Laplace-Beltrami operators with Dirichlet boundary conditions.
We introduce a {\em weak isospectral condition} (H$_1$)-(H$_2$) for the family $\Delta_t$, $t\in [0,1]$.
Generally speaking, the family $\Delta_t$, $t\in [0,1 ]$, 
is {\em weakly isospectral} if the spectrum of  $\Delta_t$ lies for all $t$ in a given  disjoint union of
infinitely many intervals $[b_k,c_k]$, $k\in\N$,  going to infinity,  polynomially separated and of
length $c_k-b_k=o(\sqrt{b_k})$ which may even tend to infinity as $k\to \infty$. 
Of course, the numbers $b_k$ and $c_k$ do not depend on $t$. 

The spectral rigidity result (see Theorem \ref{teo:rigidity-of-LBT}) obtained in Sect. \ref{sect:rigidity-of-BT} can
be formulated as follows. Consider a $C^1$ family of billiard tables $X_t$, $t\in [0,1 ]$, on the ellipsoid
such that $X_0\in {\mathcal O}(X)$ for some $X\in {\mathcal E}$,  $X_t\in {\mathcal B}$ for any $t\in [0,1]$,
and the boundaries  $\Gamma_0$ and $\Gamma_1$ are analytic. If the corresponding family of
Laplace-Beltrami operators $\Delta_t$, $t\in [0,1 ]$,  is {\em weakly isospectral}, then  $X_0=X_1$.


\section{Proof of Theorem \ref{th:hessian_elliptic_point}}\label{sect:proof_hessian}
Let $\{(x,y,p_1,p_2)\}$ be the standard coordinates on the cotangent bundle $T^*C$. Consider the Hamiltonians 
\[
H=\frac{p_1^2+p_2^2}{f(x)-q(y)},\quad F=\frac{q(y)p_1^2+f(x)p_2^2}{f(x)-q(y)}
\]
corresponding to the metric and the integral \eqref{eq:g-Liouville_form} by the Legendre transform 
\begin{equation}\label{eq:impulses}
p_1=\big(f(x)-q(y)\big){\dot x},\quad p_2=\big(f(x)-q(y)\big){\dot y}
\end{equation}
where $\dot x$ and $\dot y$ denote the components of the velocity vectors in $TC$.
For 
\[
0<h\le \max f=\alpha_0
\] 
consider the invariant with respect to the geodesic flow on $T^*C$ surface 
\[
Q_h:=\big\{H=1, F=h\big\}\subseteq T^*C.
\]
Since the variables separate one easily sees that $Q_h$ is characterized by the set of equations
\begin{equation}\label{eq:alternative}
p_1^2=f(x)-h,\quad p_2^2=h-q(y).
\end{equation}
One concludes from \eqref{eq:alternative} and the hypothesis on the functions $f$ and
$q$ that for $0<h<\alpha_0$ the surface $Q_h$ consists of two copies of $(\R/\Z)\times[-N,N]$.
The billiard reflection map at the boundary of $C$ preserves the boundary of $Q_h$ and can be used
to ``glue'' the two components of $Q_h$ into a single Liouville torus ${\widetilde Q}_h$ of the broken geodesic 
flow on $C$. By integrating the Liouville form $\kappa=p_1dx+p_2dy$ along the two cycles on 
${\widetilde Q}_h$ that are parallel to the coordinate lines on $C$, we obtain from \eqref{eq:alternative} 
that (cf. \cite{To})
\begin{equation}\label{eq:K}
L(h)=2\int_{-N}^N\sqrt{h-q(y)}\,dy
\end{equation}
and
\begin{equation}\label{eq:I}
I(h)=2\int_{x_h'}^{x_h''}\sqrt{f(x)-h}\,dx
\end{equation}
where  $0<x_h'\le 1/4\le x_h''<1/2$ are the two zeros of the equation $f(x)=h$.
Note that $x_h'=x_h''=1/4$ if and only if $h=\alpha_0$.
Since $f$ has a Morse singularity at $x_0=1/2$ there exists an orientation preserving
change of variables $p : U(0)\to V(0)$ from an open neighborhood of zero $U(0)$ 
onto an open neighborhood of zero $V(0)$ such that
\begin{equation}\label{eq:u}
x-x_0=p(u),
\end{equation}
and
\begin{equation}\label{eq:morse}
f(x)-h=(\alpha_0-h)-u^2.
\end{equation}
It follows directly from \eqref{eq:f_taylor}, \eqref{eq:u}, and \eqref{eq:morse} that
\begin{equation}\label{eq:p_Taylor}
p(u)=\frac{1}{\sqrt{-\alpha_1}}\,y+\frac{\alpha_2}{2\alpha_1^2\sqrt{-\alpha_1}}\,y^3+O(y^5).
\end{equation}
In view of \eqref{eq:I}, \eqref{eq:u}, \eqref{eq:morse}, and \eqref{eq:p_Taylor} we obtain
\begin{eqnarray*}
I(h)&=&2\int_{x_h'}^{x_h''}\sqrt{f(x)-h}\,dx\\
&=&2\int_{-\sqrt{\alpha_0-h}}^{\sqrt{\alpha_0-h}}p'(u)\sqrt{(\alpha_0-h)-y^2}\,dy\\
&=&2(\alpha_0-h)\int_{-1}^1p'\big(u\sqrt{\alpha_0-h}\big)\sqrt{1-u^2}\,du\\
&=&-\frac{\pi}{\sqrt{-\alpha_1}}(h-\alpha_0)+\frac{3\alpha_2\pi}{8\alpha_1^2\sqrt{-\alpha_1}}(h-\alpha_0)^2
+O\big((h-\alpha_0)^3\big).
\end{eqnarray*}
Hence,
\begin{equation}\label{eq:I-derivatives}
I(\alpha_0)=0,\quad\frac{dI}{dh}(\alpha_0)=-\pi/\sqrt{-\alpha_1},\quad
\frac{d^2I}{dh^2}(\alpha_0)=3\pi\alpha_2/4\alpha_1^2\sqrt{-\alpha_1}.
\end{equation}
It follows from \eqref{eq:K} and the fact that $\alpha_0>0$ that $L$ is a $C^\infty$-smooth function of $h$ in an open
neighborhood of $h=\alpha_0$. By combining this with \eqref{eq:I-derivatives} we obtain
\begin{equation}\label{eq:rotation_limit}
\frac{dL}{dI}(0)=\frac{dK}{dh}(\alpha_0)\Big\slash\frac{dI}{dh}(\alpha_0)=
-\frac{\sqrt{-\alpha_1}}{\pi}\int_{-N}^N\frac{dy}{\sqrt{\alpha_0-q(y)}}
\end{equation}
and
\begin{eqnarray}
\frac{d^2L}{dI^2}(0)\!\!&=&\!\!\left( \frac{d^2L}{dh^2}(\alpha_0)-\frac{dL}{dI}(0)\frac{d^2I}{dh^2}(\alpha_0)\right)
\Big\slash\Big(\frac{dI}{dh}(\alpha_0)\Big)^2\nonumber\\
&=&\!\!\frac{\alpha_1}{4\pi^2}\left(2\int_{-N}^N\frac{dy}{\big(\alpha_0-q(y)\big)^{3/2}}-
\frac{3\alpha_2}{\alpha_1^2}\int_{-N}^N\frac{dy}{\sqrt{\alpha_0-q(y)}}\right).\label{eq:rotation'_limit}
\end{eqnarray}
This completes the proof of Theorem \ref{th:hessian_elliptic_point}.

\section{Non-degeneracy of the Liouville billiard tables on the ellipsoid}\label{sect. LBT_on_ellipsoid}
Assume that $0<a_0<a_1<a_2$ are fixed and consider the ellipsoid in $\R^3$ with coordinates
$\{(x_0,x_1,x_2)\}$,
\begin{equation}\label{eq:ellipsoid}
\frac{x_0^2}{a_0}+\frac{x_1^2}{a_1}+\frac{x_2^2}{a_2}=1
\end{equation}
equipped with the restriction of the Euclidean metric $d\ell^2:=dx_0^2+dx_1^2+dx_2^2$.
The elliptic coordinates $\{(u_1,u_2)\}$ in the first octant are given by the formulas
\begin{eqnarray*}
x_0(u_1,u_2):=\frac{\sqrt{a_0(u_1-a_0)(u_2-a_0)}}{\sqrt{(a_1-a_0)(a_2-a_0)}}\\
x_1(u_1,u_2):=\frac{\sqrt{a_1(a_1-u_1)(u_2-a_1)}}{\sqrt{(a_1-a_0)(a_2-a_1)}}\\
x_2(u_1,u_2):=\frac{\sqrt{a_2(a_2-u_1)(a_2-u_2)}}{\sqrt{(a_2-a_0)(a_2-a_1)}}
\end{eqnarray*}
where $a_0\le u_1\le a_1$ and $a_1\le u_2\le a_2$ (see e.g. \cite[\S 3.5]{Klingenberg}).
A direct computation shows that in elliptic coordinates the metric tensor of the restriction of the
Euclidean metric $d\ell^2$ on the ellipsoid takes the form (see e.g. \cite[Proposition 3.5.4]{Klingenberg}))
\begin{equation}\label{eq:g-ellipsoid}
dg^2=(u_2-u_1)\Big(U_1(u_1)du_1^2+U_2(u_2)du_2^2\Big)
\end{equation}
where
\[
U_1(u_1):=\frac{u_1}{4(u_1-a_0)(a_1-u_1)(a_2-u_1)}
\]
and
\[
U_2(u_2):=\frac{u_2}{4(u_2-a_0)(u_2-a_1)(a_2-u_2)}\,.
\]
Consider the functions
\begin{equation}\label{eq:X}
X(u_2):=\int_{a_1}^{u_2}\frac{\sqrt{s}}{2\sqrt{(s-a_0)(s-a_1)(a_2-s)}}\,ds
\end{equation}
and
\begin{equation}\label{eq:Y}
Y(u_1):=\int_{a_1}^{u_1}\frac{\sqrt{t}}{2\sqrt{(t-a_0)(a_1-t)(a_2-t)}}\,dt\,.
\end{equation}
Note that if we pass to the new variables $x:=X(u_2)$ and $y:=Y(u_1)$,
$a_0\le u_1\le a_1$, $a_1\le u_2\le a_2$, the metric \eqref{eq:g-ellipsoid} will
take the Liouville form \eqref{eq:g-Liouville_form} on the square\footnote{Note that the variables $x$ and $y$
parametrize only the part of the ellipsoid situated in the first octant in $\R^3$.}
$0\le x\le\om_1/4$, $0\le y\le\om_2/4$ where $\om_1/4:=X(a_2)>0$ and $\om_2/4:=Y(a_0)>0$ are the periods and 
\begin{equation}\label{eq:f,g}
f(x):=X^{-1}(x)\quad\text{and}\quad q(y):=Y^{-1}(y)
\end{equation}
where $X^{-1} : [0,\om_1/4]\to[a_1,a_2]$ and $Y^{-1} : [0,\om_2/4]\to[a_0,a_1]$
are the inverses of the functions \eqref{eq:X} and \eqref{eq:Y} respectively. 
The following Remark allows us to apply Theorem \ref{th:hessian_elliptic_point} also in the case when 
the $x$-period of the covering cylinder $C$ in \eqref{eq:covering} is non-necessarily equal to one.

\begin{Rem}\label{rem:scaling}
For a given $\mu>0$ the change of variables $(x,y)\mapsto(\mu x,\mu y)$ on the cylinder 
$C\equiv(\R/\Z)\times[-N,N]$ leads to a Liouville metric $dg_\mu ^2$ on the cylinder 
$C_\mu \equiv(\R/\mu \Z)\times[-N_\mu ,N_\mu ]$, $N_\mu :=\mu  N$. 
It is easy to see that the formulas \eqref{eq:K'} and \eqref{eq:K''} continue to hold for the Liouville billiard 
corresponding to this metric.
\end{Rem}

\noindent In particular, we see that
\[
a_0\le q(y)\le a_1\quad\text{and}\quad a_1\le f(x)\le a_2
\]
where
\begin{equation}\label{eq:f,g-range}
q(\om_2/4)=a_0, q(0)=a_1\quad\text{and}\quad f(0)=a_1, f(\om_1/4)=a_2.
\end{equation}
The billiards of the first type correspond to the range $x\in[0,\om_1/4]$, $y\in[0,y_\bullet]$,
where the value $y_\bullet\in(0,\om_1/4)$ is fixed.
We will first compute the Taylor's expansion \eqref{eq:f_taylor} of $f(x)$ at $x_0=\om_1/4$ 
(see Remark \ref{rem:scaling}). In view of \eqref{eq:f,g} we obtain from \eqref{eq:X} that
\begin{equation}\label{eq:f'-relation}
\frac{df}{dx}(x)=2\frac{\sqrt{(f(x)-a_0)(f(x)-a_1)(a_2-f(x))}}{\sqrt{f(x)}},\quad x\in[0,\om_1/4].
\end{equation}
This together with \eqref{eq:f,g-range} implies that
\[
f(\om_1/4)=a_2\quad\text{and}\quad \frac{df}{dx}(\om_1/4)=0.
\]
We then differentiate the relation \eqref{eq:f'-relation} with respect to $x$ and use the chain rule
to conclude that
\[
\frac{d^2f}{dx^2}(x)=-2\,\frac{2f(x)^3-(a_0+a_1+a_2)f(x)^2+a_0a_1a_2}{f(x)^2}.
\]
By setting $x=\om_1/4$ in this formula we obtain from $f(\om_1/4)=a_2$ that
\[
\frac{d^2f}{dx^2}(\om_1/4)=-2\,\frac{(a_2-a_0)(a_2-a1)}{a_2}\,.
\]
We then continue this process to conclude that
\[
\frac{d^3f}{dx^3}(\om_1/4)=0\quad\text{and}\quad 
\frac{d^4f}{dx^4}(\om_1/4)=8\,\frac{(a_2-a_0)(a_2-a_1)(a_2^2-a_0a_1)}{a_2^3}\,.
\]
Hence,\footnote{The formulas for the derivatives were obtained in collaboration with Seth Berman as a part of his 
undergraduate research project.}
\begin{equation}\label{eq:alpha}
\left\{
\begin{array}{l}
\alpha_0=a_2,\\ 
\displaystyle \alpha_1=-\frac{(a_2-a_0)(a_2-a_1)}{a_2}<0,\\
\displaystyle\alpha_2=\frac{1}{3}\frac{(a_2-a_0)(a_2-a_1)(a_2^2-a_0a_1)}{a_2^3}>0
\end{array}
\right.
\end{equation}
and 
\begin{equation}\label{eq:kappa}
\kappa:=\frac{3\alpha_2}{\alpha_1^2}=\frac{a_2^2-a_0a_1}{a_2(a_2-a_0)(a_2-a_1)}>0
\end{equation}
where $\kappa$ is the coefficient appearing in front of the second integral in \eqref{eq:K''}.
We then pass to the new variable $t=Y^{-1}(y)$, $y\in(0,\om_2/4)$, in the integrals in \eqref{eq:K'} and \eqref{eq:K''}
to obtain that for any value of the parameter $N\in(0,\om_2/4)$,
\[
\int_{-N}^{N}\frac{dy}{(\alpha_0-q(y))^{3/2}}=\int_{\lambda}^{a_1}\frac{\sqrt{t}}{\sqrt{(t-a_0)(a_1-t)}(a_2-t)^2}\,dt
\]
and
\[
\int_{-N}^{N}\frac{dy}{\sqrt{\alpha_0-q(y)}}=\int_{\lambda}^{a_1}\frac{\sqrt{t}}{\sqrt{(t-a_0)(a_1-t)}(a_2-t)}\,dt
\]
where $\lambda\equiv Y^{-1}(N)\in(a_0,a_1)$. 
Hence,
\begin{equation}\label{eq:K'-ellipsoid}
\frac{dL}{dI}(0,\lambda)=-\frac{\sqrt{-\alpha_1}}{\pi}\int_{
\lambda}^{a_1}\frac{\sqrt{t}}{\sqrt{(t-a_0)(a_1-t)}(a_2-t)}\,dt
\end{equation}
and
\begin{equation}\label{eq:K''-ellipsoid}
\frac{d^2L}{dI^2}(0,\lambda)=\frac{\alpha_1}{4\pi^2}
\int_{\lambda}^{a_1}\frac{\sqrt{t}}{\sqrt{(t-a_0)(a_1-t)}(a_2-t)^2}\big(2-\kappa(a_2-t)\big)\,dt
\end{equation}
where $\lambda\in(a_0,a_1)$ and the coefficients $\alpha_1$ and $\kappa$ are given in \eqref{eq:alpha} and \eqref{eq:kappa}.

\begin{Rem}
Recall that the  parameter $\lambda \in(a_0,a_1)$ has a simple geometrical meaning: the corresponding Liouville billiard $X_\lambda$  on the ellipsoid is obtained by 
the intersection of the ellipsoid \eqref{eq:ellipsoid} with the interior of the hyperboloid of one sheet
\eqref{eq:confocal_family}. In this way, $\lambda \in(a_0,a_1)$ parametrizes the different billiards
of the first type on the ellipsoid \eqref{eq:ellipsoid}.
\end{Rem}

The quantity in \eqref{eq:K'-ellipsoid} when considered as a function of $\lambda \in(a_0,a_1)$ is 
strictly increasing and
\[
\frac{dL}{dI}(0,a_0+ )<\frac{dL}{dI}(0)<0
\]
where
\begin{align*}
\frac{dL}{dI}(0,a_0+)&=
-\frac{\sqrt{-\alpha_1}}{\pi}\int_{a_0}^{a_1}\frac{\sqrt{t}}{\sqrt{(t-a_0)(a_1-t)}(a_2-t)}\,dt\\
&\ge-\frac{\sqrt{-\alpha_1}}{\pi}\int_{a_0}^{a_1}\frac{\sqrt{a_1}}{\sqrt{(t-a_0)(a_1-t)}(a_2-t)}\,dt\\
&=-\frac{\sqrt{a_1}}{\sqrt{a_2}}>-1.
\end{align*}
Hence, for any $\lambda \in(a_0,a_1)$ we have that
\begin{equation}\label{eq:K'-ellipsoid_range}
0<-\frac{dL}{dI}(0,\lambda )<\sqrt{\frac{a_1}{a_2}}<1
\end{equation}
as claimed in Theorem \ref{th:non-degeneracy}. Clearly, $\lambda \mapsto -\frac{dL}{dI}(0)$ is a strictly decreasing
function of $\lambda \in(a_0,a_1)$ such that $\frac{dL}{dI}(0, a_1-)=0$.

Let us now turn our attention to the quantity in \eqref{eq:K''-ellipsoid} when considered as a function
of $\lambda\in(a_0,a_1)$. Denote by $\Psi(t)$ the integrand in \eqref{eq:K''-ellipsoid},
\begin{equation}\label{eq:Psi}
\Psi(t):=
-\frac{(-\alpha_1)}{4\pi^2}\frac{\sqrt{t}}{\sqrt{(t-a_0)(a_1-t)}(a_2-t)^2}\big(2-\kappa(a_2-t)\big)
\end{equation}
where $\alpha_1<0$.
The linear factor $2-\kappa(a_2-t)=\kappa t+(2-a_2)$ vanishes at the point
\[
t_*:=a_2-\frac{2a_2(a_2-a_0)(a_2-a1)}{a_2^2-a_0a_1}
\]
that satisfies the inequality 
\[
t_*<a_1\,.
\]
A direct computation shows that $t_*\in(a_0,a_1)$ if $a_2\in\big(a_1,a_1+\sqrt{a_1(a_1-a_0)}\big)$ and
$t_*\le a_0$ if $a_2\ge a_1+\sqrt{a_1(a_1-a_0)}$. In particular, we see form \eqref{eq:Psi} that
$\Psi(t)<0$ for $t\in(T_*,a_1)$ and $\Psi(t)>0$ for $t<T_*$. Since 
\[
\frac{d^2L}{dI^2}(0,a_1-)=0,
\]  
in order to show that 
\[
\frac{d^2L}{dI^2}(0,\lambda )<0,\quad\lambda\in(a_0,a_1),
\]
it suffices to prove that $\frac{d^2L}{dI^2}(0,a_0+)<0$. This amounts to proving that
\[
\int_{a_0}^{a_1}\frac{\sqrt{t}}{\sqrt{(t-a_0)(a_1-t)}(a_2-t)^2}\big(2-\kappa(a_2-t)\big)\,dt>0,
\]
where $\kappa>0$ is given in \eqref{eq:kappa}. The last integral equals the quantity
\begin{equation}\label{eq:E1}
\frac{2a_0\sqrt{a_1}}{a_2(a_2-a_0)(a_2-a_1)}
\left(\frac{a_2}{a_0} E\Big(\sqrt{1-\frac{a_0}{a_1}}\Big)-K\Big(\sqrt{1-\frac{a_0}{a_1}}\Big)\right)
\end{equation}
where, for $k\in(0,1)$,
\begin{equation}\label{eq:EllipticK}
K(k):=\int_0^1\frac{ds}{\sqrt{(1-s^2)(1-k^2s^2)}}
\end{equation}
is an elliptic integral of the first kind and
\[
E(k):=\int_0^1\frac{\sqrt{1-k^2s^2}}{\sqrt{1-s^2}}\,ds\,.
\]
In order to show that \eqref{eq:E1} is positive we will prove that
\begin{equation}\label{eq:E2}
\frac{a_1}{a_0}E\Big(\sqrt{1-\frac{a_0}{a_1}}\Big)-K\Big(\sqrt{1-\frac{a_0}{a_1}}\Big)>0.
\end{equation}
(Note that the coefficient $a_2/a_0$ in front of the term $E\Big(\sqrt{1-\frac{a_0}{a_1}}\Big)$ in
\eqref{eq:E1} is now replaced by the smaller constant $a_1/a_0$.)
The substitution $k:=\sqrt{1-\frac{a_0}{a_1}}$ in \eqref{eq:E2} leads to the expression
$\frac{1}{1-k^2}\,E(k)-K(k)$ that we denote by $Z(k)$, namely, 
\[
Z(k):=\frac{1}{1-k^2}\,E(k)-K(k),
\] 
and we consider it as a function of $k\in(0,1)$. We have that (see e.g. \cite[\S 22.736]{WhWo})
\[
\frac{dK}{dk}(k)=\frac{1}{k}\Big(\frac{1}{1-k^2}\,E(k)-K(k)\Big)=Z(k)/k.
\]
Hence
\[
Z(k)=k\,\frac{dK}{dk}(k)
\]
where $K(k)$ is the elliptic integral of the first kind \eqref{eq:EllipticK}.
Since
\[
\frac{dK}{dk}(k)=\int_0^1\frac{ks^2}{\sqrt{1-s^2}(1-k^2s^2)^{3/2}}\,ds>0,\quad k\in(0,1),
\]
we conclude that $Z(k)>0$ for any $k\in(0,1)$. 
Hence
\[
\frac{d^2L}{dI^2}(0,\lambda)<0
\]
for any value of $\lambda\in(a_0,a_1)$.

\medskip

Let us now consider the second type of billiards on the ellipsoid \eqref{eq:ellipsoid}.
In order to apply the formulas in Theorem \ref{th:hessian_elliptic_point} we pass
to the new variables $v_1:=-u_2$, $v_2:=-u_1$, $b_0\le v_1\le b_1$, $b_1\le v_2\le b_2$, where
\begin{equation}\label{eq:b}
b_0:=-a_2,\quad b_1:=-a_1,\quad\text{and}\quad b_2:=-a_0,
\end{equation}
and hence $b_0<b_1<b_2<0$.
Then, the metric tensor \eqref{eq:g-ellipsoid} takes the form
\begin{equation}\label{eq:g-ellipsoid2}
dg^2=(v_2-v_1)\Big(V_1(v_1)dv_1^2+V_2(v_2)dv_2^2\Big)
\end{equation}
where
\[
V_1(v_1):=\frac{(-v_1)}{4(v_1-b_0)(b_1-v_1)(b_2-v_1)}
\]
and
\[
V_2(v_2):=\frac{(-v_2)}{4(v_2-b_0)(v_2-b_1)(b_2-v_2)}\,.
\]
We then set
\begin{equation}\label{eq:X2}
X(v_2):=\int_{b_1}^{v_2}\frac{\sqrt{-s}}{2\sqrt{(s-b_0)(s-b_1)(b_2-s)}}\,dt,
\end{equation}
\begin{equation}\label{eq:Y2}
Y(v_2):=\int_{v_1}^{b_1}\frac{\sqrt{-s}}{2\sqrt{(s-b_0)(b_1-s)(b_2-s)}}\,dt\,.
\end{equation}
and pass to the variables $x:=X(v_2)$ and $y:=Y(v_1)$ where $b_0\le v_1\le b_1$,
$b_1\le v_2\le b_2$. As in the case of the billiards of the first type, 
the metric tensor \eqref{eq:g-ellipsoid2} in these coordinates has a Liouville form \eqref{eq:g-Liouville_form} 
where $f(x):=X^{-1}(x)$, $q(y):=Y^{-1}(y)$, and $0\le x\le\om_1/4$, $0\le y\le\om_2/4$, with
$\om_1/4:=X(b_2)$ and $\om_2/4:=Y(b_0)$. In particular,
\[
b_0\le q(y)\le b_1\quad\text{and}\quad b_1\le f(x)\le b_2
\]
and
\begin{equation}\label{eq:f,g-range2}
q(\om_2/4)=b_0, q(0)=b_1\quad\text{and}\quad f(0)=b_1, f(\om_1/4)=b_2.
\end{equation}
In order to compute the Taylor's expansion \eqref{eq:f_taylor} of $f(x)$ at $x_0=\om_1/4$ we proceed
in the same way as in the case of the billiards of the first type: From \eqref{eq:X2} and the fact that
$f(x)=X^{-1}(x)$ we conclude that
\[
\frac{df}{dx}(x)=\frac{2\sqrt{(f(x)-b_0)(f(x)-b_1)(b_2-f(x))}}{\sqrt{-f(x)}},\quad x\in[0,\om_1/4).
\]
We then differentiate this expression with respect to $x$ and use the chain rule and the fact that
$f(\om_1/4)=b_2$ (cf. \eqref{eq:f,g-range2}) to conclude that
\begin{equation}\label{eq:alpha2}
\left\{
\begin{array}{l}
\alpha_0=b_2=-a_0<0,\\ [0.3cm]
 \alpha_1=\frac{(b_2-b_0)(b_2-b_1)}{b_2}=-\frac{(a_2-a_0)(a_1-a_0)}{a_0}<0,\\ [0.3cm]
 \alpha_2=-\frac{1}{3}\frac{(b_2-b_0)(b_2-b_1)(b_0b_1-b_2^2)}{b_2^3}
=\frac{1}{3}\frac{(a_2-a_0)(a_1-a_0)(a_2a_1-a_0^2)}{a_0^3}>0
\end{array}
\right.
\end{equation}
and 
\begin{equation}\label{eq:kappa2}
\kappa:=\frac{3\alpha_2}{\alpha_1^2}=\frac{a_2a_1-a_0^2}{a_0(a_2-a_0)(a_1-a_0)}>0\,.
\end{equation}
By passing to the variable $t=Y^{-1}(y)$, $y\in(0,\om_2/4)$, in the integrals appearing in \eqref{eq:K'} and \eqref{eq:K''}
we obtain that for any value of the parameter $N\in(0,\om_2/4)$,
\begin{align*}
\int_{-N}^{N}\frac{dy}{(\alpha_0-q(y))^{3/2}}&=\int_{(-\lambda)}^{b_1}\frac{\sqrt{-t}}{\sqrt{(t-b_0)(b_1-t)}(b_2-t)^2}\,dt\\
&=\int_{a_1}^T\frac{\sqrt{s}}{\sqrt{(a_2-s)(s-a_1)}(s-a_0)^2}\,ds
\end{align*}
and
\begin{align*}
\int_{-N}^{N}\frac{dy}{\sqrt{\alpha_0-q(y)}}&=\int_{(-\lambda)}^{b_1}\frac{\sqrt{-t}}{\sqrt{(t-b_0)(b_1-t)}(b_2-t)}\,dt\\
&=\int_{a_1}^T\frac{\sqrt{s}}{\sqrt{(a_2-s)(s-a_1)}(s-a_0)}\,ds
\end{align*}
where $\lambda\equiv -Y^{-1}(N)\in(a_1,a_2)$. 
This and Theorem \ref{th:hessian_elliptic_point} then imply that for any $\lambda\in(a_1,a_2)$ we have that
\begin{equation}\label{eq:K'-ellipsoid2}
\frac{dL}{dI}(0,\lambda)=-\frac{\sqrt{-\alpha_1}}{\pi}\int_{a_1}^\lambda\frac{\sqrt{s}}{\sqrt{(a_2-s)(s-a_1)}(s-a_0)}\,ds
\end{equation}
and
\begin{equation}\label{eq:K''-ellipsoid2}
\frac{d^2L}{dI^2}(0,\lambda)=\frac{\alpha_1}{4\pi^2}
\int_{a_1}^\lambda\frac{\sqrt{s}}{\sqrt{(a_2-s)(s-a_1)}(s-a_0)^2}\big(2-\kappa(s-a_0)\big)\,ds
\end{equation}
where $\alpha_1$ and $\kappa$ are given in \eqref{eq:alpha2} and \eqref{eq:kappa2}.

\begin{Rem}
As in the case of the billiards of the first type, the parameter $\lambda\in(a_1,a_2)$ has the following 
geometrical meaning: the corresponding Liouville billiard $X_\lambda$ on the ellipsoid is the intersection of 
the ellipsoid \eqref{eq:ellipsoid} with the hyperboloid of two sheets \eqref{eq:confocal_family}.
\end{Rem}

The integrand in \eqref{eq:K''-ellipsoid2} equals
\[
\Phi(s):=\frac{(-\alpha_1)}{4\pi^2}
\frac{\sqrt{s}}{\sqrt{(a_2-s)(s-a_1)}(s-a_0)^2}\big(\kappa(s-a_0)-2\big)
\]
where $\alpha_1<0$.
The linear factor $\kappa(s-a_0)-2=\kappa s-(\kappa a_0+2)$ vanishes at the  point
\[
s_*:=a_0+\frac{2 a_0(a_2-a_0)(a_1-a_0)}{a_2 a_1-a_0^2}
\]
A direct computation shows that
\[
s_*<a_1\,.
\]
This implies that $\Phi(s)>0$ for $s\in(a_1,a_2)$. Hence,
\[
\frac{d^2L}{dI^2}(0,\lambda)>0
\]
for any value of $\lambda\in(a_1,a_2)$.

It follows from \eqref{eq:K'-ellipsoid2}, the equality 
$\int_{a_1}^{a_2}\frac{ds}{\sqrt{(a_2-s)(s-a_1)}}=\pi$, and the fact that 
the integration in \eqref{eq:K'-ellipsoid2} is over the interval $a_1\le s\le \lambda$ with $\lambda<a_2$,
that for any $\lambda\in(a_1,a_2)$ we have
\[
0<-\frac{dL}{dI}(0,\lambda)<\sqrt{\frac{a_2(a_2-a_0)}{a_0(a_1-a_0)}}
\]
where $-\frac{dL}{dI}(0,\lambda)$ is a strictly increasing function of $\lambda\in(a_1,a_2)$ that satisfies
$\frac{dL}{dI}(0,a_1+)=0$ and
\[
\sqrt{\frac{a_1(a_1-a_0)}{a_0(a_2-a_0)}}<-\frac{dL}{dI}(0,a_2-).
\]
This completes the proof of Theorem \ref{th:non-degeneracy} for the billiards on the ellipsoid.

\section{Spectral rigidity of billiard tables}\label{sect:rigidity-of-BT}
We first recall from \cite{PT2}  a spectral rigidity result for $C^1$  deformations of a  billiard table $X$  with boundary
$\Gamma:=\partial X$ in an ambient Riemannian manifold $(\widetilde X,g$) of dimension two equipped with a metric $g$.
To this end we suppose that $X$ has an elliptic  $4$-elementary  bouncing ball geodesic $\gamma_1$, or equivalently that
the symplectic map $P:=B^2$ admits a Birkhoff normal form at any of  its vertices  $p_\ast\in B^\ast\Gamma$ given by
\begin{equation}\label{eq:BNF}
 (\theta,I)\mapsto\Big(\theta+\frac{dL}{dI}(I) + O(|I|^{3/2}),I+O(|I|^{5/2})\Big)
\end{equation}
where $(\theta,I)\in (\R/\Z)\times \R_+$ are suitable polar symplectic coordinates, $I(p_\ast)=0$, 
$L$ is a smooth function such that  $L(0)=0$,  $\frac{dL}{dI}(0)\neq 0 \ \textrm{mod}\, 1$, 
and it satisfies  the non-resonance condition \eqref{eq:4-elementary}. The Birkhoff normal form \eqref{eq:BNF} is said to
be {\em non-degenerate} if
\begin{equation}\label{eq:non-degenerate-BNF}
\frac{d^2L}{dI^2}(0)\neq 0.
\end{equation}
Let $X$ be a compact billiard table in an ambient Riemannian manifold $(\widetilde X,g)$, $\mbox{dim}\, \widetilde X=2$, 
with boundary $\Gamma:=\partial X$. 
\begin{Def}\label{def:C-1-families}
We say that 
the family of billiard tables $X_t$, $t\in [0,1]$,  is a $C^1$ deformation of  $X$ in $\widetilde X$ if  
\begin{itemize}
\item[$\rhd$] $X_t$ are smooth compact connected submanifolds of $\widetilde X $ of dimension two with
boundaries $\Gamma_t=\partial X_t$  given by a $C^1$ family of  embeddings $\psi_t:\Gamma\to \widetilde X$,
$t\in [0,1]$, such that $\psi_0={\rm id}$,  $\Gamma_t=\psi_t(\Gamma)$, and the map
$[0,1]\ni t\mapsto \psi_t\in C^\infty(\Gamma, \widetilde X) $ is $C^1$-smooth; 
\item[$\rhd$]  for any $t_0\in [0,1]$ and $\varrho_0$ in the interior of  $X_{t_0}$ there exists a neighborhood
$U$ of $t_0$ in  $[0,1]$ such that $\varrho_0\in X_t$ whenever $t\in U$. 
\end{itemize}
\end{Def}
If  $X_{t_0}$  is given then the second condition determines uniquely  the interior of the billiard table  $X_t$ from its
boundary $\Gamma_t=\partial X_t$ for any $t\in [0,1]$. 

Denote by $\Delta_t$ the ``geometric'' Laplace-Beltrami operator corresponding to the Riemannian manifold
$(X_t,g)$ with Dirichlet  boundary conditions.  This is a self-adjoint operator in $L^2(X_t)$ with discrete spectrum
accumulating at $+\infty$.  We define a  {\em weak isospectral condition} as follows. 

Consider a union ${\mathcal U}$ of infinitely many disjoint intervals $[b_k,c_k]$ going to infinity, of
length $o(\sqrt{b_k})$, and polynomially separated. 

More precisely, fix  two positive constants $d\ge 0$ and $A>0$, and suppose that
\begin{enumerate}
\item[$(\mbox{H}_1)$] {\em  ${\mathcal U}\subseteq(0,\infty) $ 
is a union of infinitely many disjoint  intervals $[b_k,c_k]$, $k\in\N$,  such that 
\begin{enumerate}
\item[$\rhd$] $\displaystyle{\lim\,  b_k\, =\, \lim\,  c_k\,  =\,  +\infty}$;
\item[$\rhd$] $\displaystyle\lim\,   \frac{c_k-b_k}{\sqrt{b_k}}\,  =\,  0$; 
\item[$\rhd$]
$b_{k+1} - c_{k}\, \ge \, A c_k^{-d}$ \quad for any $k\in \N$. 
\end{enumerate}}
\end{enumerate}
Given a set ${\mathcal U}$ satisfying $(\mbox{H}_1)$, 
we impose the following {\em ``weak isospectral assumption''} 
\begin{enumerate}
\item[$(\mbox{H}_2)$] {\em There exists   $a\ge 1$ such that  \quad 
${\rm Spec}\left(\Delta_t \right)\, \cap [a,+\infty)\  
\subseteq \   {\mathcal U}\quad  \forall\,  t \in [0,1]   \, .$}
\end{enumerate}
Note that the length of the intervals $[b_k,c_k]$ can increase and even go  to infinity as $k\to \infty$ but not faster than
$o(\sqrt{b_k})$.  Using the asymptotic behavior  of the eigenvalues $\lambda_j$ as $j\to \infty$ one can  show that
(H$_1$)-(H$_2$) are ``natural'' for any $d> n/2$ and $c>0$. By ``natural'' we mean that  the usual isospectral condition 
\[
{\rm Spec}\left(\Delta_t \right)\, = \, 
{\rm Spec}\left(\Delta_0 \right) \quad \forall\,  t \in [0,1] 
\]
implies that there exists $a\ge 1$ and  a family of infinitely many disjoint intervals $[b_k,c_k]$ such that (H$_1$)-(H$_2$) are 
satisfied -- see \cite{PT3}, Lemma 2.2, for details. 

Let us define a class of billiard tables as follows. Let $(\widetilde X,g)$, $\mbox{dim}\, \widetilde X=2$, be a Riemannian
manifold of dimension two. Suppose that there exist two geodesics $\tilde\gamma_k$, $k=1,2$, 
intersecting transversely at a point $\varrho_0$ and two commuting  involutive isometries 
${\mathcal J}_k : {\widetilde X}\to{\widetilde X}$, $k=1,2$, such that the set of fixed points of ${\mathcal J}_k$
coincides with $\tilde\gamma_k$. 
Consider the family ${\mathcal B}$ of compact billiard tables $X$ in $(\widetilde X,g)$ containing $\varrho_0$ such that
the boundary $\Gamma=\partial X$ of $X$ is connected, it intersects each $\tilde\gamma_k$, $k=1,2$  at exactly
two points and $\Gamma$ is invariant with respect to ${\mathcal J}_k$, $k=1,2$. Then the set of fixed points of 
${\mathcal J}_k$, $k=1,2$, in $X$ defines a bouncing ball geodesic of $(X,g)$  that we denote by $\gamma_k$. 

The following result has been obtained in  \cite{PT2}. 

\begin{Th}\label{theo:elliptic} 
Let $X_0\in{\mathcal B}$. Assume that the  broken geodesic $\gamma_1$ given by the set of fixed points of ${\mathcal J}_1$
in $X_0$ is elliptic $4$-elementary and that  the corresponding Poincare map $P:=B^2$ admits a non-degenerate
Birkhoff normal form. Suppose that $X_t\in{\mathcal B}$, $t\in [0,1]$,  is a $C^1$ deformation of $X_0$ satisfying
the weak isospectral condition $(\mbox{H}_1)-(\mbox{H}_2)$. Then   $\gamma_1$ is a bouncing ball geodesic of
$(X_t,g)$ for each $t\in [0,1]$ and $\Gamma_t$ has a contact of infinite order to $\Gamma_0$ at the vertexes of
$\gamma_1$. In particular, $X_0=X_1$ if both $\Gamma_0$ and $\Gamma_1$ are
analytic. 
\end{Th}
The proof of Theorem \ref{theo:elliptic} is based on a Kolmogorov-Arnold-Moser (KAM) theorem and a quasi-mode construction 
for $C^1$ families of billiard tables (see \cite{PT2}).

\subsection{Spectral rigidity of billiard tables on the ellipsoid}\label{sec:spectral_rigidity}
We are going to apply Theorem \ref{theo:elliptic} to the Liouville billiard tabels $X^\lambda$ of  first and second type on the
ellipsoid $\widetilde X:=E(a_0,a_1,a_2)$ . Consider the involutions ${\mathcal J}^{\alpha}_k$, $k=1,2$, $\alpha={\rm I, I\!I}$,
of the ellipsoid  $E(a_0,a_1,a_2)$ defined by
\begin{equation}\label{eq:involutions_ellipsoid}
\begin{array}{lcr}
{\mathcal J}^{\rm I}_1(x_0,x_1,x_2)=(x_0,x_1,-x_2) , \ {\mathcal J}^{\rm I}_2(x_0,x_1,x_2)=(x_0,-x_1,x_2)\\ [0.3cm]
{\mathcal J}^{\rm I\!I}_1(x_0,x_1,x_2)=(-x_0,x_1,x_2) , \ {\mathcal J}^{\rm I\!I}_2(x_0,x_1,x_2)=(x_0,-x_1,x_2)
\end{array}
\end{equation}
The involutions ${\mathcal J}^{\rm I}_k$,  $k=1,2$, coincide with the involutions ${\mathcal J}_k$  of the billiard tables
$X^\lambda$ of {\em  first type }given by   Remark  \ref{rem:involutions}, while  ${\mathcal J}^{\rm I\!I}_k$,  $k=1,2$,
are the corresponding involutions for billiards of  {\em second type} on the ellipsoid. We have $\varrho_0= \pm (\sqrt{a_0},0,0)$
when $\alpha={\rm I}$ and $\varrho_0= \pm (0,0, \sqrt{a_2})$ when $\alpha={\rm I\!I}$. 
Let us  denote by ${\mathcal B}^\alpha$,  $\alpha= {\rm I,I\!I}$,   the corresponding classes of billiard tables on
the ellipsoid $E(a_0,a_1,a_2)$.  Denote by  $\mathcal{E}^{\rm I}$ and  $\mathcal{E}^{\rm I\!I}$ the families
of Liouville billiard tables  $X^\lambda$  of first   and second type  on  $E(a_0,a_1,a_2)$  for which the non-resonance
condition \eqref{eq:4-elementary} is satisfied. Let $\Lambda^{\rm I}\subseteq(a_0,a_1)$ and
$\Lambda^{\rm I\!I}\subseteq (a_1,a_2)$ be the sets of the corresponding  values of $\lambda$. 
Thus we get the parametrization 
\[
\Lambda^\alpha\ni \lambda \mapsto X^\lambda\in \mathcal{E}^{\alpha},\ \alpha={\rm I,I\!I}.
\]
Corollary \ref{coro:non-degeneracy} says that 
$\Lambda^{\rm I}=(a_0,a_1)$ if $\sqrt{\frac{a_1}{a_2}}<1/4$. In general, the complement of $\Lambda ^{\rm I}$ in $(a_0,a_1)$
consists of at most five exceptional points. These exceptional values of $\lambda$ are obtained by means of  \eqref{eq:K'-ellipsoid}. 
The complement of $\Lambda^{\rm I\!I}$ in $(a_1,a_2)$ consists of at most finitely many exceptional points.   
These exceptional values of $\lambda$ are obtained by means of  \eqref{eq:K'-ellipsoid2}. 

The non-resonance condition \eqref{eq:4-elementary} and the non-degeneracy condition
\eqref{eq:non-degenerate-BNF} of the Birkhoff normal form \eqref{eq:BNF} are open  in the $C^5$ topology. 
Then for any $X^\lambda\in \mathcal{E}^\alpha$, $\alpha= {\rm I,I\!I}$, there exists a  $\delta$-neighborhood
${\mathcal O}^\alpha (\lambda)\subseteq {\mathcal B}^\alpha$ of 
$X^\lambda$ in the $C^5$ topology such that for each $X\in {\mathcal O}^\alpha (\lambda)$ the corresponding  geodesic
$\gamma_1$ is elliptic $4$-elementary and the Birkhoff normal form is non-degenerate. The set ${\mathcal O}^\alpha (\lambda)$ is
defined as follows. We say that a billiard table $ X\in {\mathcal B}^\alpha$  belongs to 
${\mathcal O}^\alpha (\lambda)$ if there exists an embedding $f:\Gamma^\lambda \to E(a_0,a_1,a_2)$,
such that $f(\Gamma^\lambda)=\Gamma:=\partial X$, and 
\[
\|f-{\rm id}\|_{C^5(\Gamma^\lambda)}\le\delta
\]
where $\delta>0$ is chosen sufficiently small. 
Applying Theorem \ref{theo:elliptic} we obtain the following  spectral rigidity result. 
\begin{Th}\label{teo:rigidity-of-LBT}
Let $\lambda\in \Lambda^\alpha$, $\alpha= {\rm I,I\!I}$, and $X\in \mathcal{O}^\alpha(\lambda)$. Denote by $\gamma_1$
the corresponding elliptic geodesic. 
Suppose that $X_t\in{\mathcal B}$, $t\in [0,1]$,  is a $C^1$ deformation of $X_0:=X$ satisfying
the weak isospectral condition $(\mbox{H}_1)-(\mbox{H}_2)$. Then   $\gamma_1$ is a bouncing ball geodesic of
$X_t$ for each $t\in [0,1]$ and $\Gamma_t:=\partial X_t$ has a contact of infinite order to $\Gamma_0$ at the vertexes of
$\gamma_1$. In particular, $X_0=X_1$ if  both $\Gamma_0$ and $\Gamma_1$ are
analytic. 
\end{Th}

\subsection{Preceding results}\label{sect:preceding_results}
The inverse spectral problem for   billiard tables  was one of the main research topics  of Steve Zelditch. 
There is a long list of excellent results of Steve  \cite{H-Z1}-\cite{H-Z5} and \cite{Z1}-\cite{Z6}  (some of them in collaboration with Hamid Hezari) about the spectral 
determination and the  spectral rigidity of billiard tables, where fine microlocal analysis and dynamical systems are 
applied to investigating  the singularities of the wave-trace. It was proved in \cite{Z4} and \cite{Z5} that a generic real analytic plain domain with
one up-down symmetry is uniquely determined by the spectrum of the corresponding Laplace operator. 
A similar result for billiard tables in $\R^n$, $n\ge 2$, with a  $(\Z/2\Z)^n$ group of symmetries was
obtained in  \cite{H-Z2}.  Under some generic non-degeneracy assumptions, it is proved in \cite{H-Z4} that real analytic centrally symmetric plane domains are determined by their Dirichlet (resp. Neumann) spectra in this class. 
Hezari and Zelditch \cite{H-Z2} proved that ellipses are infinitesimally spectrally
rigid in the class of smooth billiard tables  with left-right and up-down symmetries. 
In \cite{H-Z3} Hezari and Zelditch prove that ellipses with small eccentricity are uniquely determined by
the Dirichlet (or Neumann) Laplace spectrum. Moreover, they obtain spectral rigidity in the class of
$\Z/2\Z$ - axially symmetric smooth domains close in $C^8$ to the unit disc. Closely related to these results
are also the articles \cite{ASK}, \cite{DKW}, \cite{KS}, on the Birkhoff conjecture and the dynamical rigidity.

Infinitesimal spectral rigidity of Kolmogorov non-degenerate Liouville billiard tables in dimensions two and
three (in particular for the ellipse and the   ellipsoid) has been obtained in \cite{PT2} under the
{\em weak isospectral condition}. The proof is based on a KAM theorem for $C^1$ families of billiard tables,
a quasi-mode construction, and the injectivity of the Radon transform  investigated in  \cite{PT4}.
Spectral rigidity results of Liouville billiard tables on surfaces of constant curvature similar to that in
Theorem \ref{teo:rigidity-of-LBT} are obtained in \cite{PT2}.

\end{document}